\documentclass[a4paper]{amsart}
\usepackage{amssymb,amsthm,amsmath,amscd,amsfonts,bbm,mathrsfs}

\usepackage[cmtip,all]{xy}

\theoremstyle{plain}
\newtheorem{Lemma}{Lemma}
\newtheorem{Thm}[Lemma]{Theorem}
\newtheorem*{Thm*}{Theorem}
\newtheorem{Prop}[Lemma]{Proposition}
\newtheorem{Cor}[Lemma]{Corollary}
\theoremstyle{definition}
\newtheorem{Def}[Lemma]{Definition}
\theoremstyle{remark}
\newtheorem*{Remark*}{Remark}
\newtheorem*{Remarks*}{Remarks}

\numberwithin{Lemma}{section}
\numberwithin{equation}{section}

\newcommand{\CCC}{\mathcal{C}}
\newcommand{\EEE}{\mathcal{E}}
\newcommand{\FFF}{\mathcal{F}}
\newcommand{\GGG}{\mathscr{G}}

\newcommand{\NNN}{\mathcal{N}}
\newcommand{\OOO}{\mathcal{O}}
\newcommand{\PPP}{\mathscr{P}}
\newcommand{\ZZZ}{\mathscr{Z}}
\newcommand{\Fa}{\mathfrak{a}}
\newcommand{\Fb}{\mathfrak{b}}
\newcommand{\Fm}{\mathfrak{m}}
\newcommand{\Gm}{{\mathbb{G}_m}}
\newcommand{\hGa}{{\widehat{\mathbb{G}}_a}}
\newcommand{\hGm}{{\widehat{\mathbb{G}}_m}}
\newcommand{\II}{{\mathbb{I}}}
\newcommand{\QQ}{{\mathbb{Q}}}
\newcommand{\WW}{{\mathbb{W}}}
\newcommand{\ZZ}{{\mathbb{Z}}}

\newcommand{\hatP}{\widehat P}
\newcommand{\hatQ}{\widehat Q}
\newcommand{\hatV}{\widehat V}
\newcommand{\hatW}{\widehat W}
\newcommand{\hatZ}{\widehat Z}

\newcommand{\tCCC}{\widetilde\CCC}
\newcommand{\tEEE}{\widetilde\EEE}
\newcommand{\tL}{\widetilde L}
\newcommand{\tP}{\widetilde P}
\newcommand{\tQ}{\widetilde Q}
\newcommand{\tT}{\widetilde T}

\newcommand{\fep}{\mbox{$\frac 1p$}}

\DeclareMathOperator{\Bil}{Bil}
\DeclareMathOperator{\Ext}{Ext}

\DeclareMathOperator{\Ker}{Ker}
\DeclareMathOperator{\Lie}{Lie}
\DeclareMathOperator{\Gal}{Gal}
\DeclareMathOperator{\Hom}{Hom}

\DeclareMathOperator{\Spec}{Spec}
\DeclareMathOperator{\red}{red}
\DeclareMathOperator{\can}{can}
\DeclareMathOperator{\hex}{hex}
\DeclareMathOperator{\id}{id}
\DeclareMathOperator{\et}{et}
\DeclareMathOperator{\op}{op}
\DeclareMathOperator{\nil}{nil}
\DeclareMathOperator{\mult}{mult}

\DeclareMathOperator{\incl}{incl}
\DeclareMathOperator{\Spf}{Spf}
\DeclareMathOperator{\BT}{BT}

\newcommand{\EExt}{\underline\Ext}
\newcommand{\HHom}{\underline\Hom}

\begin{document}

\title{A Duality Theorem for Dieudonn\'{e} Displays}
\author{Eike Lau}
\date{\today}
\address{Fakult\"{a}t f\"{u}r Mathematik,
Universit\"{a}t Bielefeld, D-33501 Bielefeld, Germany}
\email{lau@math.uni-bielefeld.de}

\begin{abstract}
We show that the Zink equivalence 
between $p$-divisible groups and Dieudonn\'{e} displays 
over a complete local ring with perfect residue field
of characteristic $p$ is compatible with duality. 
The proof relies on a new explicit formula for the 
$p$-divisible group associated to a Dieudonn\'{e} display.
\end{abstract}

\maketitle

\section*{Introduction}

Let $R$ be a complete local ring with maximal ideal $\Fm$ 
and perfect residue field $k$ of positive characteristic $p$.
If $p=2$ we assume that $pR=0$.

As a generalisation of classical Dieudonn\'e\ theory,
Th.~Zink defines in \cite{Zink-DDisp} a category of Dieudonn\'e\
displays over $R$ and shows that it is equivalent 
to the category of $p$-divisible groups over $R$. 
In the present article we give a unified formula for the 
group associated to a Dieudonn\'e\ display and 
apply it to show that the equivalence is compatible 
with the natural duality operations on both sides.
This is not clear from the original construction
because that depends on decomposing a $p$-divisible group
into its etale and infinitesimal part,
which is not preserved under duality.

Let us recall the definition of a Dieudonn\'e\ display.
There is a unique subring $\WW(R)$ of the Witt ring $W(R)$
that is stable under its Frobenius $f$ and Verschiebung $v$,
that surjects onto $W(k)$, and that contains an element 
$x\in W(\Fm)$ if and only if the components of $x$ converge
to zero $\Fm$-adically.
In \cite{Zink-DDisp} the ring $\WW(R)$ is denoted $\hatW(R)$.
Let $\II_R$ be the kernel of the natural homomorphism 
$\WW(R)\to R$. 
A Dieudonn\'e\ display over $R$ is a quadruple 
$$
\PPP=(P,Q,F,F_1)
$$ 
where $P$ is a finite free $\WW(R)$-module,
$Q$ a submodule containing $\II_RP$ such that 
$P/Q$ is a free $R$-module, 
$F:P\to P$ and $F_1:Q\to P$ are $f$-linear maps such 
that $F_1(v(w)x)=wF(x)$ for $x\in P$ and $w\in\WW(R)$,
and the image of $F_1$ generates $P$.
Dieudonn\'e\ displays over $k$ are equivalent to Dieudonn\'e\
modules $(P,F,V)$ where $Q=V(P)$ and $F_1=V^{-1}$.

Our formula is based on viewing both $p$-divisible
groups and the modules $P$, $Q$ as abelian sheaves 
for the flat topology on the opposite
category of all $R$-algebras $S$
with the following properties: 
the nilradical $\NNN(S)$ is a nilpotent ideal, 
it contains $\Fm S$, and $S/\NNN(S)$ 
is a union of finite dimensional $k$-algebras;
see section \ref{Se1} for details.
With that convention, the equivalence functor $\BT$ from 
Dieudonn\'e\ displays to $p$-divisible groups is given by
\begin{equation}
\label{EqI0}
\tag{$\star$}
\BT(\PPP)=[Q\xrightarrow{F_1-\incl}P]\overset{L}{\otimes}\QQ_p/\ZZ_p
\end{equation}
where $[Q\to P]$ is a complex of sheaves in degrees $0,1$.
In other words, the cohomology of the right hand side of \eqref{EqI0} 
vanishes outside degree zero and the zeroth cohomology is 
the $p$-divisible group associated to $\PPP$. 
Instead of the flat topology one could also use the
ind-etale topology, but for some arguments the former 
is more convenient.

Before stating the main result let us
recall the duality of Dieudonn\'e\ displays.
We have the special Dieudonn\'e\ display 
$\GGG_m=(\WW(R),\II_R,f,v^{-1})$ that
corresponds to the group $\hGm$.
A bilinear form $\PPP'\times\PPP\to\GGG_m$ 
is a bilinear map $\alpha:P'\times P\to\WW(R)$
satisfying $\alpha(x',x)=v(\alpha(F_1'x',F_1x))$ 
for $x'\in Q'$ and $x\in Q$. 
For every $\PPP$ there is a dual $\PPP^t$
equipped with a perfect bilinear form 
$\PPP^t\times\PPP\to\GGG_m$,
which determines $\PPP^t$ uniquely.
The Serre dual of a $p$-divisible group $G$
is denoted $G^\vee$.

\begin{Thm*}
For every Dieudonn\'e\ display $\PPP$ over $R$ 
there is a natural isomorphism 
$$
\Psi:\BT(\PPP^t)\cong\BT(\PPP)^\vee.
$$
\end{Thm*} 

The proof is independent of the fact that 
the functor $\BT$ from Dieudonn\'e\ displays to 
$p$-divisible groups defined by \eqref{EqI0} is 
actually an equivalence.
Let us indicate how to define the homomorphism $\Psi$.
Denote by $Z(\PPP)$ the complex $[Q\to P]$ in \eqref{EqI0}.
To the tautological bilinear form $\PPP^t\times\PPP\to\GGG_m$
one can directly assign a homomorphism of complexes 
$Z(\PPP^t)\otimes Z(\PPP)\to Z(\GGG_m)$, 
which gives after tensoring twice with 
$\QQ_p/\ZZ_p$ a homomorphism
$$
\BT(\PPP^t)\overset L\otimes\BT(\PPP)\to
\BT(\GGG_m)\overset L\otimes\QQ_p/\ZZ_p\cong\hGm[1].
$$
By the cohomological theory of biextensions,
such a homomorphism is equivalent to 
a homomorphism $\Psi$ as above.
That $\Psi$ is an isomorphism must be shown only if the group 
$\BT(\PPP)$ is etale or of multiplicative type or bi-infinitesimal.
The first two cases are straightforward; 
the bi-infinitesimal case relies on the theorem 
of Cartier \cite{Cartier} 
on the Cartier dual of the Witt ring functor. 

Over arbitrary rings in which $p$ is nilpotent, 
infinitesimal $p$-divisible groups are equivalent to {\em displays} 
according to \cite{Zink-Disp} and \cite{Lau}. 
The bi-infinitesimal case of the above theorem is closely related 
to the duality theorem in \cite{Zink-Disp} for the display 
associated to a bi-infinitesimal $p$-divisible group. 
This in turn has been anticipated by Norman \cite{Norman}
who shows a similar duality theorem for the Cartier module 
of a bi-infinitesimal $p$-divisible group, provided the module 
is displayed (which is always the case by the said equivalence).
These duality results all depend on the theory of biextensions
developed in \cite{Mumford}, that appears here
in the cohomological form it was given in SGA 7.

The present proof that the functor $\BT$ defined by \eqref{EqI0} 
is an equivalence of categories consists in verifying that 
it reproduces the equivalence constructed in \cite{Zink-DDisp}.
However, it should be possible to relate the crystals associated 
to a Dieudonn\'e\ display $\PPP$ and to the $p$-divisible group 
$\BT(\PPP)$. 
Then the fact that $\BT$ is an equivalence will follow
directly from the Grothendieck-Messing deformation theory of
$p$-divisible groups \cite{Messing-Crys},
and the duality theorem for 
Diedonn\'{e} displays will be related
to the crystalline duality theorem \cite{BBM}. 
We hope to return to this point later.
Let us also note that Caruso \cite{Caruso}
proved a duality theorem for Breuil modules of finite
flat $p$-group schemes \cite{Breuil}
by using the crystalline duality theorem.
Breuil modules of $p$-divisible
groups are related to Dieudonn\'{e} displays 
by \cite{Zink-Windows}.

This text is organised as follows.
In section \ref{Se1} the formula for $\BT$ 
is explained, in section \ref{Se2}
it is shown to give an equivalence of categories,
in section \ref{Se3} the duality theorem is proved,
and section \ref{Se4} is concerned with 
functoriality in the base.
In an appendix we discuss briefly the deformational
duality theorem \cite{Mazur-Messing} since variants
of it are used in the text.

\medskip
I am grateful to W.~Messing and Th.~Zink 
for many valuable discussions and to one 
of the referees for pointing out that section
\ref{Se4} was missing.

\section{Exposition of the main formula }
\label{Se1}

We begin with a number of general definitions and notations.
Let $p$ be a prime.
For any ring $A$ let $W(A)$ be the ring of $p$-Witt vectors and
$I_A$ the kernel of the first Witt polynomial $w_0:W(A)\to A$.
If $A$ is perfect of characteristic $p$, $I_A$ is generated by $p$.
Let $f$ be the Frobenius of $W(A)$ and $v$ the Verschiebung. 
If $\Fa\subset A$ is a nilpotent ideal, 
let $\hatW(\Fa)\subseteq W(\Fa)$ be the subgroup of 
Witt vectors with only finitely many non-zero components. 
More generally, if $A$ is $\Fa$-adically complete and separated,
let $\hatW(\Fa)\subseteq W(\Fa)$ be the subgroup of Witt vectors
whose components converge to zero $\Fa$-adically;
in other words, $\hatW(\Fa)=\varprojlim \hatW(\Fa/\Fa^n)$. 
In any case $\hatW(\Fa)$ is an ideal in $W(A)$.

\begin{Def}
\label{DefAdm}
Let $A$ be a ring and $\Fa\subset A$ an ideal.
The pair $(A,\Fa)$ is called admissible if $A$ is
$\Fa$-adically complete and separated and $A/\Fa$ 
is perfect of characteristic $p$.
If $p=2$ we also require that $pA=0$.
\end{Def}

\begin{Lemma}
If $(A,\Fa)$ is admissible
then there is a unique $f$-stable subring $\WW(A)$ of $W(A)$
such that $\WW(A)\cap W(\Fa)=\hatW(\Fa)$ and 
$\WW(A)$ maps surjectively onto $W(A/\Fa)$. 
The subring $\WW(A)$ is also stable under $v$.
\end{Lemma}

This is proved in \cite{Zink-DDisp} if $A$ is noetherian and $A/\Fa$ 
is a field, but neither of these assumptions is used in the proof.
$\WW(A)$ is constructed as follows: Since $A/\Fa$ is perfect,
the projection $W(A)\to W(A/\Fa)$ has a unique splitting,
necessarily $f$-equivariant, thus an $f$-equivariant
decomposition of abelian groups $W(A)\cong W(A/\Fa)\oplus W(\Fa)$,
under which $\WW(A)$ is mapped to $W(A/\Fa)\oplus\hatW(\Fa)$.
The condition that $2$ is invertible or zero in $A$ is only
needed to guarantee that $\WW(A)$ is $v$-stable.
We have 
$$
\WW(A)=\raisebox{0em}[0em][0em]{$\varprojlim\WW(A/\Fa^n)$}
$$ 
by uniqueness or by the construction. 
Let $\II_A$ be the kernel of $w_0:\WW(A)\to A$.

\begin{Def}
Assume that $(A,\Fa)$ is admissible.
A Dieudonn\'{e} display over $A$ is a quadruple
$\PPP=(P,Q,F,F_1)$ such that
\begin{quote}
$P$ is a finitely generated projective $\WW(A)$-module, \\
$Q$ is a submodule of $P$ containing $\II_AP$, \\
$P/Q$ is projective as an $A$-module, \\
$F:P\to P$ and $F_1:Q\to P$ are $f$-linear maps,\\
$F_1(v(w)x)=wF(x)$ for $w\in\WW(A)$ and $x\in P$,\\
$F_1(Q)$ generates $P$ as a $\WW(A)$-module.
\end{quote}
These axioms also imply $F(x)=pF_1(x)$ for $x\in Q$. 
\end{Def}

\begin{Remarks*}
(1) {}
Every pair of $\WW(A)$-modules $(P,Q)$ satisfying the 
first three of the above conditions admits a decomposition 
$P=L\oplus T$ such that $Q=L\oplus \II_AT$, 
called normal decomposition. 
Its existence is straightforward if $\Fa=0$, thus $A$ perfect;
if $\Fa$ is nilpotent one can use that $\hatW(\Fa)$ is
nilpotent as well; the general case follows by passing 
to the projective limit.

(2) {}
If $F:M\to N$ is an $f$-linear homomorphism of $\WW(A)$-modules,
let $M^{(1)}=\WW(A)\otimes_{f,\WW(A)}M$, and let 
$F^\sharp:M^{(1)}\to N$ be the linearisation of $F$. 
In analogy with \cite{Zink-Disp} Lemma 9,
the structure of a Dieudonn\'e\ display on a pair $(P,Q)$ 
as above with given normal decomposition $P=L\oplus T$
is equivalent to the isomorphism 
$$
(F_1^\sharp,F^\sharp):L^{(1)}\oplus T^{(1)}\xrightarrow{\sim} P.
$$

(3) {}
We have the following notion of base change.
If $(A,\Fa)$ and $(B,\Fb)$ are admissible, a ring homomorphism
$g:A\to B$ with $g(\Fa)\subseteq\Fb$ induces a ring homomorphism
$\WW(g):\WW(A)\to\WW(B)$.
The base change of a Dieudonn\'e\ display $\PPP$ over $A$ by $g$ 
is then $\PPP_B=(P_B,Q_B,F_B,F_{1,B})$ where
$$
P_B=\WW(B)\otimes_{\WW(A)}P,\qquad
Q_B=\Ker(P_B\to B\otimes_AP/Q),
$$
and $F_B$, $F_{1,B}$ are the unique 
$f$-linear extensions of $F$, $F_1$,
whose existence follows from a normal decomposition 
as explained in \cite{Zink-Disp} Definition 20.

(4) {}
For a Dieudonn\'e\ display $\PPP$ over $A$ 
there is a unique $\WW(A)$-linear map 
$V^\sharp:P\to P^{(1)}$ 
such that $V^\sharp(F_1x)=1\otimes x$ for $x\in Q$,
cf.\ \cite{Zink-Disp} Lemma 10. Uniqueness is clear;
if $P=L\oplus T$ is a normal decomposition, 
$V^\sharp$ can be defined to be  
$$
P\xrightarrow{(F_1^\sharp,F^\sharp)^{-1}}
L^{(1)}\oplus T^{(1)}\xrightarrow{(1,p)}
L^{(1)}\oplus T^{(1)}=P^{(1)}.
$$
We have $F^\sharp V^\sharp=p$ and $V^\sharp F^\sharp=p$.
If $A$ is perfect, $F_1$ is bijective, and its inverse
defines an $f^{-1}$-linear map $V:P\to P$
whose linearisation is $V^\sharp$.
\end{Remarks*}

Assume now that $R$ is a local ring with maximal ideal $\Fm$ 
and residue field $k$ such that $(R,\Fm)$ is admissible,
i.e.\ $R$ is $\Fm$-adically complete, $k$ is perfect
of characteristic $p$, and $p=2$ implies $pR=0$. 

\begin{Def} 
Let $\CCC_R$ be the category of all $R$-algebras $S$
such that the nilradical $\NNN(S)$ is nilpotent,
$\NNN(S)$ contains $\Fm S$, 
and $S_{\red}=S/\NNN(S)$ is a union of finite dimensional, 
necessarily etale, $k$-algebras.
\end{Def}

The last condition implies that $S_{\red}$ is perfect, 
hence $(S,\NNN(S))$ is admissible, and $\WW(S)$ is defined. 
The following stability properties of $\CCC_R$ are easily established:
If $S'\leftarrow S\to S''$ are morphisms in $\CCC_R$, then
$S'\otimes_SS''$ lies in $\CCC_R$; if $S\in\CCC_R$ and $S\to S'$
is a finite ring homomorphism, then $S'\in\CCC_R$; if $S\in\CCC_R$
and $S\to S_1\to S_2\to\ldots$ is an infinite sequence of etale
ring homomorphisms, then $\varinjlim S_i$ lies in $\CCC_R$.

Let $\raisebox{0em}[0em][0em]{$\tCCC_R$}$ 
be the category of abelian sheaves on $\CCC_R^{\op}$
for the flat topology, i.e.\ coverings are all
faithfully flat ring homomorphisms in $\CCC_R$. 
The category of $p$-divisible groups over $R$ is 
naturally a full exact subcategory of $\tCCC_R$ 
that is stable under extensions.
If $\PPP$ is a Dieudonn\'e\ display over $R$,
base change of Dieudonn\'e\ displays makes $P$ and $Q$ into 
abelian presheaves on $\CCC_R^{\op}$, i.e.\ for $S\in\CCC_R$
we put $P(S)=P_S$ and $Q(S)=Q_S$.
Note that the presheaf $Q$ is determined by the
modules $Q\subseteq P$ but not by the module
$Q$ alone. The homomorphisms $F$ and $F_1$
induce homomorphisms of the associated presheaves
which we denote by the same letters.
By the following lemma, $P$ and $Q$ are in fact sheaves.

\begin{Lemma}
For a faithfully flat homomorphism $S\to T$ in $\CCC_R$ the
natural sequence
$0\to \WW(S) \to \WW(T) \rightrightarrows \WW(T\otimes_ST)$
is exact.
\end{Lemma}

\begin{proof}
The analogous assertion with $W$ in place of $\WW$ is clear,
cf.\ \cite{Zink-Disp} Lemma 30.
One easily checks that $\WW(S)=\WW(T)\cap W(S)$, 
and the lemma follows.
\end{proof}

\begin{Def}
\label{DefBT}
If $\PPP$ is a Dieudonn\'e\ display over $R$, let 
$$
Z(\PPP)=[Q\xrightarrow{F_1-\incl}P]
$$
as a complex in $\tCCC_R$ in degrees $0,1$ and 
$$
\BT(\PPP)=Z(\PPP)\otimes^L\QQ_p/\ZZ_p
$$
in the derived category $D(\tCCC_R)$.
\end{Def}

Explicitly $\BT(\PPP)$ can be represented 
by the tensor product of complexes
$Z(\PPP)\otimes[\ZZ\to\ZZ[\frac{1}{p}]]$ 
sitting in degrees $-1,0,1$.

\begin{Thm}
\label{Thm1}
Suppose $R$ is an admissible local ring.
For every Dieudonn\'e\ display $\PPP$ over $R$, 
$\BT(\PPP)$ is a $p$-divisible group, 
i.e.\ $H^i(\BT(\PPP))$ vanishes for $i\neq 0$ 
and is a $p$-divisible group for $i=0$. 
The functor $\BT$ induces an equivalence of 
exact categories 
$$
\{\text{Dieudonn\'e\ displays over $R$}\}\cong
\{\text{$p$-divisible groups over $R$}\}
$$
that coincides with the equivalence in \cite{Zink-DDisp}. 
The height of $\BT(\PPP)$ is equal to the rank of $P$,
and there is a natural isomorphism 
$\Lie(\BT(\PPP))\cong P/Q$.
\end{Thm}

Here the additive category of Dieudonn\'e\ displays 
is made into an exact category by declaring a 
short sequence to be exact if it is exact 
on the $P$'s and on the $Q$'s.
Let us stress again that the only new aspect in 
Theorem \ref{Thm1} is the formula for the functor $\BT$.
It will be proved in the next section.

\begin{Remark*}
The functor $\BT$ is also compatible
with base change, see section \ref{Se4}.
\end{Remark*}

\section{Proof of the main formula}
\label{Se2}

Let $R$, $\Fm$, $k$ be as before.
We begin with recalling some definitions and results from 
\cite{Zink-DDisp} that are stated there only if $R$ is artinian, 
but if $\Fm$ is nilpotent the arguments apply without change,
and the general case follows by passing to the limit
since Dieudonn\'e\ displays over $R$ are equivalent
to compatible systems of Dieudonn\'e\ displays 
over $R/\Fm^n$ for $n\geq 1$.

A Dieudonn\'e\ display $\PPP$ over $R$ is called 
{\em etale} if $V^\sharp$ is an isomorphism, 
{\em of multiplicative type} if $F^\sharp$ is an isomorphism,
and {\em $V$-nilpotent} or {\em $F$-nilpotent} if $V^\sharp$ or $F^\sharp$ 
is topologically nilpotent for the adic topology on $\WW(R)$ 
defined by the ideal $\hatW(\Fm)+\II_R$.
$\PPP$ is etale if and only if $Q=P$ 
and of multiplicative type if and only if $Q=\II_RP$,
see \cite{Zink-DDisp} Definitions 13 \& 14.
Etale or multiplicative or $V$-nilpotent 
or $F$-nilpotent Dieudonn\'e\ displays over $R$ 
are equivalent to compatible systems of the same objects
over $R/\Fm^n$ for $n\geq 1$ because each of these conditions 
holds for $\PPP$ over $R$ if and only if it holds 
for the base change $\PPP_k$ over $k$.

By \cite{Zink-DDisp} Propositions 15, 16 \& 17, there are no 
non-trivial homomorphisms between etale and $V$-nilpotent
or between multiplicative and $F$-nil\-potent Dieudonn\'e\ displays 
in either direction, moreover for every $\PPP$ there are 
unique and functorial exact sequences of Dieudonn\'e\ displays 
\begin{gather}
\label{EqS1}
0\to\PPP^{V-{\nil}}\to\PPP\to\PPP^{\et}\to 0 \\ 
\label{EqS2}
0\to\PPP^{\mult}\to\PPP\to\PPP^{F-{\nil}}\to 0
\end{gather}
such that $\PPP^{\et}$, $\PPP^{\mult}$, $\PPP^{V-{\nil}}$,
$\PPP^{F-{\nil}}$ are of the designated types.
The corresponding assertions for $p$-divisible groups are well-known:
Let us call a $p$-divisible group $G$ over $R$ infinitesimal if
$G(\bar k)=(0)$, 
i.e.\ if $G$ is infinitesimal as a group over $\Spf R$.
Then for every $G$ there is a unique and functorial exact sequence
of $p$-divisible groups 
\begin{equation}
\label{EqE}
0\to G^{\inf}\to G\to G^{\et}\to 0 \\
\end{equation}
such that $G^{\et}$ is etale and $G^{\inf}$ infinitesimal,
moreover by rigidity there are no non-zero homomorphisms between 
etale and infinitesimal $p$-divisible groups over $R$ in 
either direction.

The equivalence between Dieudonn\'e\ displays and 
$p$-divisible groups in \cite{Zink-DDisp} is obtained by showing
that $V$-nilpotent or etale Dieudonn\'e\ displays are equivalent 
to infinitesimal or etale $p$-divisible groups, respectively,
and by providing an explicit isomorphism 
$\Ext^1(\PPP,\PPP')\cong\Ext^1(G,G')$ 
if $\PPP$ is an etale and $\PPP'$ a $V$-nilpotent 
Dieudonn\'e\ display and $G$, $G'$ are the associated groups.
In order to prove Theorem \ref{Thm1} we show that the functor
$\BT$ reproduces the given equivalences in the etale 
and $V$-nilpotent case and that 
it induces the given isomorphism on $\Ext^1$.

As a preparation we define for every Dieudonn\'e\ display $\PPP$ 
over $R$ an exact sequence of complexes in $\tCCC_R$ of the 
following type.
\begin{equation}
\label{EqD}
0\to Z_\NNN(\PPP)\to Z(\PPP)\to\bar Z(\PPP)\to 0
\end{equation}
For $S\in\CCC_R$ let $\bar P(S)=P(S_{\red})$ and 
$\bar Q(S)=Q(S_{\red})$. Then $\bar P$ and $\bar Q$ are sheaves
on $\CCC_R^{\op}$ because for every faithfully flat 
ring homomorphism $S\to T$ in $\CCC_R$ 
the induced homomorphism $S_{\red}\to T_{\red}$ 
is also faithfully flat (all $S_{\red}$-modules are flat),
and we have
$T_{\red}\otimes_{S_{\red}}T_{\red}=(T\otimes_ST)_{\red}$. 
Let  
$$
\bar Z(\PPP)=[\bar Q\xrightarrow{F_1-\incl}\bar P],
$$
and let $Z_\NNN(\PPP)=[Q_\NNN\xrightarrow{F_1-\incl}P_\NNN]$
be the kernel of $Z(\PPP)\to\bar Z(\PPP)$, explicitly
$$
P_\NNN(S)=\hatW(\NNN(S))\otimes_{\WW(R)}P,\qquad
Q_\NNN=P_\NNN\cap Q.
$$
Note that \eqref{EqD} is already exact on the level of presheaves.

\subsection{The infinitesimal case}

\numberwithin{Lemma}{subsection}

\begin{Prop}
\label{Pr2}
Theorem \ref{Thm1} holds for $V$-nilpotent Dieudonn\'e\ displays and 
infinitesimal $p$-divisible groups.
If $\PPP$ is $V$-nilpotent and $G$ is the associated group, 
we have a natural quasi-isomorphism $Z(\PPP)\simeq G[-1]$. 
\end{Prop}

Let us recall the proof of the infinitesimal case 
in \cite{Zink-DDisp}.
To a Dieu\-don\-n\'e\ display $\PPP$ over $R$ one  
associates a display $\FFF\!\PPP=(P',Q',F',F_1')$, 
where $P'=W(R)\otimes_{\WW(R)}P$, and 
$Q'$ is the kernel of the natural map $P'\to P/Q$.
The functor $\FFF$ induces an equivalence between 
$V$-nilpotent Dieudonn\'e\ displays and 
$V$-nilpotent displays. 
This is tautological if $R=k$;  
if the maximal ideal $\Fm$ is nilpotent, 
the assertion follows by deformations;
the general case by a limit argument.

On the other hand, if $\PPP=(P,Q,F,F_1)$ 
is a {\em display} over $R$ and $N$ a nilpotent 
(non-unitary) $R/\Fm^n$-algebra for some $n$, let us write:
\begin{align*}
\hatP(N)&{}=\hatW(N)\otimes_{W(R)}P \\
\hatQ(N)&{}=\Ker(\hatP(N)\to N\otimes_R P/Q) \\
\hatZ(\PPP,N)&{}=[\hatQ(N)\xrightarrow{F_1-\incl}\hatP(N)]
\end{align*}
Then by \cite{Zink-Disp} Theorem 81 \& Corollary 89,
the group $H^0\hatZ(\PPP,N)$ vanishes, and 
the functor $N\mapsto H^1\hatZ(\PPP,N)$ is represented
by a formal group over $\Spf R$ that is $p$-divis\-ible if $\PPP$
is $V$-nilpotent. By op.~cit.\ \S 3.3,
the functor $\PPP\mapsto H^1\hatZ(\PPP,\underline{\hspace{1.5ex}}\,)$
induces an equivalence between $V$-nilpotent displays 
and infinitesimal $p$-divisible groups over $R$;
more precisely, Corollary 95 and a limit argument reduce 
this to the case $R=k$, which is covered by Proposition 102.

If $\PPP$ is a display over $R$, let $\ZZZ(\PPP)$ 
denote the complex in $\tCCC_R$ given by 
$\ZZZ(\PPP)(S)=\hatZ(\PPP,\NNN(S))$ for $S\in\CCC_R$.
For every Dieudonn\'e\ display $\PPP$ over $R$ we have 
an obvious isomorphism 
$$
Z_\NNN(\PPP)\cong\ZZZ(\FFF\PPP).
$$
It follows that $H^0 Z_\NNN(\PPP)$ vanishes, and the functor
$H^1Z_\NNN$ defines an equivalence between $V$-nilpotent 
Dieudonn\'e\ displays and infinitesimal $p$-divisible groups.
This is the equivalence of \cite{Zink-DDisp}.
Here $H^1Z_{\NNN}(\PPP)$ in the sense of presheaves or sheaves
is the same, i.e.\ the presheaf $H^1$ is already a sheaf.

\begin{Lemma}
\label{Le2}
If $\PPP$ is a $V$-nilpotent Dieudonn\'e\ display, 
then $\bar Z(\PPP)$ is acyclic.
\end{Lemma}

\begin{proof}
For $S\in\tCCC_R$ the complex 
$[F_1-\incl:Q_{S_{\red}}\to P_{S_{\red}}]$ 
is isomorphic to 
$[\id-V:P_{S_{\red}}\to P_{S_{\red}}]$ where $V=F_1^{-1}$. 
Since $V$ is topologically nilpotent, $\id-V$ is bijective.
\end{proof}

For every $K\in D(\tCCC_R)$ 
the obvious homomorphisms of complexes 
\begin{equation}
\label{Eq2}
\QQ_p/\ZZ_p\xleftarrow{\;\simeq\;}
[\ZZ\to\ZZ[\fep]]\longrightarrow\ZZ[1]
\end{equation}
(where $\simeq$ means quasi-isomorphism)
induce a morphism $\pi_K:K\otimes^L\QQ_p/\ZZ_p\to K[1]$.
It is an isomorphism if all local sections of $H^*K$ 
are annihilated by powers of $p$.

\begin{proof}[Proof of Proposition \ref{Pr2}]
By the above discussion and Lemma \ref{Le2} we have an equivalence 
$\PPP\mapsto G$ between $V$-nilpotent Dieudonn\'e\ displays
and infinitesimal $p$-divisible groups such that
$Z(\PPP)\simeq Z_\NNN(\PPP)\simeq G[-1]$.
The isomorphism $\pi_{G[-1]}$ then 
gives $\BT(\PPP)\cong G$.
Finally $\Lie(G)\cong P/Q$ by \cite{Zink-Disp} (158).
\end{proof}

\begin{Remark*}
The Dieudonn\'e\ display associated to $\hGm=\Gm[p^\infty]$ 
is given by
\begin{equation}
\label{EqG}
\GGG_m=(\WW(R),\II_R,f,v^{-1})
\end{equation}
where $v^{-1}$ is the inverse of the bijective homomorphism
$v:\WW(R)\to\II_R$.
In fact, for every nilpotent $R/\Fm^n$-algebra
$N$ there is an exact sequence
\begin{equation}
\label{Eqhex}
0\to\hatW(N)\xrightarrow{1-v}\hatW(N)
\xrightarrow{\hex}\hGm(N)\to 0
\end{equation}
where $\hex$ is given by the Artin-Hasse exponential 
evaluated at $t=1$, see \cite{Zink-Disp} p.\ 108, 
hence $\BT(\GGG_m)\cong H^1Z_\NNN(\GGG_m)\cong\hGm$.
\end{Remark*}

\subsection{The etale case}

If $G$ is an etale $p$-divisible group over $R$,
let $T_pG=\varprojlim G[p^n]$ in $\tCCC_R$.
The obvious sequences $0\to T_pG\to T_pG\to G[p^n]\to 0$
are exact as arbitrary ind-etale coverings exist in $\CCC_R$
and give an isomorphism $T_pG\otimes^L\QQ_p/\ZZ_p\cong G$.

\begin{Prop}
\label{Pr3}
Theorem \ref{Thm1} holds in the etale case.
If $G$ is the etale $p$-divisible group 
associated to an etale Dieudonn\'e\ display $\PPP$,
we have a natural quasi-iso\-morphism $T_pG\simeq Z(\PPP)$.
\end{Prop}

Note that an etale Dieudonn\'e\ display is the same as a pair
$(P,F_1)$, where $P$ is a finitely generated free
$\WW(R)$-module and $F_1^\sharp:P^{(1)}\to P$ an isomorphism;
we have $Q=P$ and $F=pF_1$.
The complex $Z(\PPP)$ takes the form $[F_1-\id:P\to P]$.

Again we have to recall the equivalence $\PPP\mapsto G$
between etale Dieudon\-n\'e\ displays and etale 
$p$-divisible groups from \cite{Zink-DDisp}. 
By op.\ cit.\ Theorem 5,
etale Dieudonn\'e\ displays over $R$ and over $k$ are equivalent.
The analogous assertion for etale $p$-divisible groups
and for truncated etale $p$-divisible groups is well-known.
Hence it suffices to define the equivalence $\PPP\mapsto G$
over $k$; there it is given by the isomorphism of
$\Gal(\bar k/k)$-modules
\begin{equation}
\label{EqT}
T_pG(\bar k)=(W(\bar k)\otimes_{W(k)}P)^{F_1=\id}.
\end{equation}

Let us reformulate this a little.
Over every ring $A$ of characteristic $p$,
truncated etale $p$-divisible groups of level $n$
are equivalent to pairs $(M,\Phi)$, where $M$
is a  finitely generated projective $W_n(A)$-module and
$\Phi:M^{(1)}\to M$ is an isomorphism. 
The group associated to $(M,\Phi)$ 
is the sheaf that maps an $A$-algebra $B$
to $(M\otimes_AB)^{\Phi=\id}$.

Assume now that $\PPP\mapsto G$ as above and that $pR=0$.
Then the pair $(M,\Phi)$ associated to the truncated
etale group $G[p^n]$ is equal to
$(W_n(R)\otimes_{\WW(R)}P,F_1)$.
In fact, to prove this one may pass to $k$,
where the assertion follows from \eqref{EqT}.
For $S\in\CCC_R$ we get
\begin{equation}
\label{EqT2}
G[p^n](S)=(W_n(S)\otimes_{\WW(R)}P)^{F_1=\id}.
\end{equation}

\begin{Lemma}
\label{Le3}
If $\PPP$ is an etale Dieudonn\'e\ display, then $Z_\NNN(\PPP)$ is acyclic.
\end{Lemma}

\begin{proof}
Since $F_1:P\to P$ is $f$-linear, the induced map
$F_1:P_\NNN\to P_\NNN$ is elementwise nilpotent,
so $F_1-\id:P_\NNN\to P_\NNN$ is bijective.
\end{proof}

\begin{proof}[Proof of Proposition \ref{Pr3}]
If $\PPP$ is an etale Dieudonn\'e\ display over $R$ 
and $G$ is the associated etale $p$-divisible group, 
the sheaves $H^i\bar Z(\PPP)$ are computed as follows.
For $S\in\CCC_R$ we have, using \eqref{EqT2},
\begin{align*}
H^0\bar Z(\PPP)(S)&{}= 
\varprojlim(W_n(S_{\red})\otimes_{\WW(R)}P)^{F_1=\id} \\
&{}=\varprojlim G[p^n](S_{\red})=
\varprojlim G[p^n](S)=T_p G(S).
\end{align*}
The sheaf $H^1\bar Z(\PPP)$ vanishes because 
every element of $P(S_{\red})$ has an inverse image under
$F_1-\id$ after passing to an ind-etale covering of
$S_{\red}$, 
which lifts to an ind-etale covering of $S$.
By Lemma \ref{Le3} we obtain
$Z(\PPP)\simeq\bar Z(\PPP)\simeq T_pG$, hence 
$\BT(\PPP)\cong G$.
\end{proof}

\subsection{Calculation of extensions} 

For an etale $p$-divisible group $H$ over $R$
there is a natural exact sequence in $\tCCC_R$ 
\begin{equation}
\label{Eq3}
0\to T_pH\to T_pH\otimes\ZZ[\fep]\to H\to 0.
\end{equation}
If $G$ is another $p$-divisible group over $R$, 
we obtain a connecting homomorphism:
\begin{equation*}
\delta:\Hom(T_pH,G)\to\Ext^1(H,G)
\end{equation*}

\begin{Remark*}
We have $\Hom(T_pH,G)\cong\varinjlim\Hom(H[p^n],G)$
because $T_pH$ is representable in $\CCC_R$, so every
homomorphism $T_pH\to G$ factors over some $G[p^n]$,
hence also over $H[p^n]$.
Using that isomorphism, $\delta$ can be defined 
by the obvious exact sequences 
$0\to H[p^n]\to H\to H\to 0$ instead of \eqref{Eq3}.
\end{Remark*}

\begin{Prop}
If $H$ is etale and $G$ is infinitesimal then $\delta$ is bijective.
\end{Prop}

This is \cite{Zink-DDisp} Proposition 19.
See Proposition \ref{PrX} for a more general statement.

We turn to extensions of Dieudonn\'e\ displays. 
Let $\PPP^{\et}$ be an etale and $\PPP^{\nil}$ a $V$-nilpotent
Dieudonn\'e\ display over $R$ and let $H=\BT(\PPP^{\et})$ and
$G=\BT(\PPP^{\nil})$ be the associated $p$-divisible groups.
For every extension of Dieudonn\'e\ displays
\begin{equation}
\label{Eq11}
0\to\PPP^{\nil}\to\PPP\to\PPP^{\et}\to 0
\end{equation}
the resulting exact sequence of complexes
\begin{equation}
\label{Eq12}
0\to Z(\PPP^{\nil})\to Z(\PPP)\to Z(\PPP^{\et})\to 0
\end{equation} 
gives rise to a connecting homomorphism
\begin{equation}
\label{Eq13}
T_pH\cong H^0Z(\PPP^{\et})\to H^1Z(\PPP^{\nil})\cong G
\end{equation}
where the outer isomorphisms are provided by
Propositions \ref{Pr2} \& \ref{Pr3}. 
This construction defines a homomorphism  
$$
\gamma:\Ext^1(\PPP^{\et},\PPP^{\nil})\to\Hom(T_pH,G).
$$

\begin{Prop}
The homomorphism $\gamma$ is bijective.
\end{Prop}

This is a reformulation of \cite{Zink-DDisp} Proposition 18
and its proof. 
The equivalence between Dieudonn\'e\ displays and $p$-divisible
groups in op.\ cit.\ is defined by the composite isomorphism
$\delta\gamma:\Ext^1(\PPP^{\et},\PPP^{\nil})\to\Ext^1(H,G)$. 

\begin{proof}[Proof of Theorem \ref{Thm1}]
By Propositions \ref{Pr2} \& \ref{Pr3} we know that 
$\BT(\PPP)$ is a $p$-divisible group of the correct 
height if $\PPP$ is etale or $V$-nilpotent. 
The same is true for general $\PPP$  
because the exact sequence \eqref{EqS1}
gives rise to a distinguished triangle 
$\BT(\PPP^{V-{\nil}})\to\BT(\PPP)\to\BT(\PPP^{\et})\to^+$.
Similarly we see that $\BT$ preserves arbitrary exact sequences.

The main point to be shown is that for 
$H=\BT(\PPP^{\et})$ and $G=\BT(\PPP^{\nil})$ 
as above, the homomorphism 
$\Ext^1(\PPP^{\et},\PPP^{\nil})\to\Ext^1(H,G)$
induced by $\BT$ coincides with the isomorphism $\delta\gamma$.

Let us first describe the action of $\BT$ on $\Ext^1$.
Let $T=T_pH$.
If an extension \eqref{Eq11} is given, 
then using the quasi-isomorphisms 
$Z(\PPP^{\nil})\simeq G[-1]$ and $Z(\PPP^{\et})\simeq T$,
the corresponding extension \eqref{Eq12}
determines the following distinguished triangle in $D(\tCCC_R)$,
where $g$ is the negative of the connecting homomorphism \eqref{Eq13}.
\begin{equation}
\label{EqX}
G[-1]\to Z(\PPP)\to T\xrightarrow{g} G
\end{equation}
The functor $\BT$ applied to \eqref{Eq11} 
results in \eqref{EqX}${}\otimes^L\QQ_p/\ZZ_p$,
which under the identification 
$\pi_{G[-1]}:G[-1]\otimes^L\QQ_p/\ZZ_p\cong G$ 
induced by \eqref{Eq2} takes the following form.
\begin{equation}
\label{EqY}
G\to\BT(\PPP)\to H\xrightarrow{g'}G[1]
\end{equation}
Here $g'$ is the composition
$$
H\cong T\otimes^L\QQ_p/\ZZ_p\xrightarrow{g\otimes\id}
G\otimes^L\QQ_p/\ZZ_p\xrightarrow{-\pi_G}G[1]
$$
because $-\pi_{G}$ gets identified with $\pi_{G[-1]}[1]$ 
under the natural isomorphism 
$G\otimes^L\QQ_p/\ZZ_p\cong G[-1]\otimes^L\QQ_p/\ZZ_p[1]$,
the sign arising from the transposition automorphism of 
$\ZZ[1]\otimes\ZZ[1]$ which is $-\id$.

Consider now $\delta\gamma$.
Since the extension \eqref{Eq3}
corresponds to the distinguished triangle 
$$
T\to T\otimes\ZZ[\fep]\to T\otimes\QQ_p/\ZZ_p
\xrightarrow{\pi_T}T[1]
$$
and since $g$ is the negative of \eqref{Eq13},
the image of \eqref{Eq11} under $\delta\gamma$ 
corresponds to a triangle of the form \eqref{EqY}
with $g'=-g[1]\pi_T$. 
It remains to show that $g[1]\pi_T=\pi_G(g\otimes\id)$,
which is clear. 

Finally, the inverse functor of $\BT$ preserves exact sequences
because this is true over $k$, and a short sequence of 
Dieudonn\'e\ displays over $R$ is exact if and only if 
it is exact over $k$ by Nakayama's lemma.
The existence of a natural isomorphism $\Lie(\BT(\PPP))\cong P/Q$
follows from the $V$-nilpotent case since
$\Lie(\BT(\PPP))=\Lie(\BT(\PPP^{\nil}))$ and $P/Q=P^{\nil}/Q^{\nil}$.
\end{proof}

\numberwithin{Lemma}{section}
\setcounter{Lemma}{\value{subsection}}

A $V$-nilpotent and $F$-nilpotent Dieudonn\'e\ display is called bi-nilpotent.

\begin{Cor}
\label{Cor1}
The group $\BT(\PPP)$ is etale or of multiplicative type or bi-infinitesimal
if and only if $\PPP$ is etale or of multiplicative type or bi-nilpotent,
respectively. 
\end{Cor}

\begin{proof}
Since a $p$-divisible group over $R$ is etale or of multiplicative type 
if and only if its dimension is equal to zero or its height,  
the etale and multiplicative case of the corollary are a direct
consequence of the last assertion of Theorem \ref{Thm1}.
The third case follows because bi-infinitesimal groups and 
bi-nilpotent Dieudonn\'e\ displays are both characterised by having
neither non-trivial sub-objects of multiplicative type nor etale quotients. 
\end{proof}

\begin{Remark*}
The proof of Corollary \ref{Cor1} requires only that
the functor $\BT$ is an equivalence over $k$, which 
is classical, because the relevant properties of 
Dieudonn\'e\ displays and $p$-divisible groups 
over $R$ depend only on their fibres over $k$.
\end{Remark*}

We conclude this section with a remark on 
the topology chosen on $\CCC_R^{\op}$.

First let us note that $p$-divisible groups over $R$ form a full 
exact subcategory of the category of sheaves on $\CCC_R^{\op}$
for the etale topology, i.e.\ every short exact sequence 
$0\to H\to E\to G\to 0$ of $p$-divisible groups
gives rise to an exact sequence of etale sheaves. 
Indeed, recall that according to 
\cite{Mazur-Messing} Lemma 10.12, 
every $H$-torsor is formally smooth.
Thus the sequence admits a set-theoretical section if $G$ is
infinitesimal, or if $G$ is etale and $H$ is infinitesimal 
since then it splits over $k$. 
By using the functorial decomposition \eqref{EqE}
if follows that the sequence of ind-etale sheaves is exact.
However, multiplication by $p$ is a surjective 
endomorphism of the etale sheaf given by a 
$p$-divisible group only if the group is etale.

For a Dieudonn\'e\ display $\PPP$ over $R$ let  
$$
Y(\PPP)=[Q\xrightarrow{F_1-1}P]
\otimes[\ZZ\to\ZZ[\fep]]
$$
as a complex of presheaves on $\CCC_R^{\op}$
concentrated in degrees $-1$, $0$, $1$. 
This is in fact a complex of flat sheaves as 
$P(S)\otimes\ZZ[\frac 1p]=Q(S)\otimes\ZZ[\frac 1p]
=W(S_{\red})[\frac 1p]\otimes_{\WW(R)}P$
for $S\in\CCC_R$, which clearly is a sheaf.
The following refinement of Theorem \ref{Thm1}
is an immediate consequence of its proof.

\begin{Cor}
\label{CorTop}
The ind-etale cohomology sheaf $H^i(Y(\PPP))$
is naturally isomorphic to $\BT(\PPP)$ for $i=0$
and vanishes for $i\neq 0$. 
\end{Cor}

\section{Duality}
\label{Se3}

Let $R$, $\Fm$, $k$ be as before.
In order to define the dual of a Dieudonn\'e\ display 
we need the following notion of bilinear forms.
Recall the definition of $\GGG_m$ in \eqref{EqG}.
If $\PPP$ and $\PPP'$ are Dieudonn\'e\ displays over $R$,
a bilinear form $\alpha:\PPP\times\PPP'\to\GGG_m$ is by definition
a $\WW(R)$-bilinear map $\alpha:P\times P'\to \WW(R)$ such that 
\begin{align}
\label{EqB1}
v(\alpha(F_1x,F_1'x')) & {}=\alpha(x,x') \\
\intertext{for $x\in Q$ and $x'\in Q'$.
This also implies for
$y\in P$, $y'\in P'$, and $x,x'$ as before:}
\label{EqB2}
\alpha(F_1x,F' y') & ={}f(\alpha(x,y')) \\
\label{EqB3}
\alpha(F y,F_1'x') & {}=f(\alpha(y,x')) \\
\label{EqB4}
\alpha(F y,F' y') & {}=p f(\alpha(y,y'))
\end{align}
Let $\Bil(\PPP\times\PPP',\GGG_m)$ denote 
the abelian group of bilinear forms.

\begin{Def}
For every Dieudonn\'e\ display $\PPP$ the contravariant functor 
$\PPP'\mapsto\Bil(\PPP'\times\PPP,\GGG_m)$ is represented
by a Dieudonn\'e\ display $\PPP^t$, called the dual of $\PPP$.
\end{Def}

This is analogous to the case of displays,
\cite{Zink-Disp} Definition 19.
Let us make the definition of $\PPP^t$ more explicit,
which can also be used to show that $\PPP^t$ exists. 
For a $\WW(R)$-module $M$ let $M^\vee$ be the 
module of linear maps $M\to\WW(R)$. Then 
$$
\PPP^t=(P^\vee,\tQ,F',F_1')
$$ 
where $\tQ=\{x\in P^\vee\mid x(Q)\subseteq\II_R\}$.
If $P=L\oplus T$ is a normal decomposition, i.e.\
$Q=L\oplus\II_RT$, accordingly $P^\vee=L^\vee\oplus T^\vee$ 
and $\tQ=\II_RL^\vee\oplus T^\vee$, the formulas
\eqref{EqB1} to \eqref{EqB4} determine that 
the composition
$$
(L^\vee)^{(1)}\oplus (T^\vee)^{(1)}
\xrightarrow{(F',F_1')^\sharp} P^\vee
\xrightarrow{(F_1^\sharp,F^\sharp)^\vee}
(L^{(1)})^\vee\oplus (T^{(1)})^\vee
$$
is the tautological isomorphism, i.e.\ when passing from
$\PPP$ to $\PPP^t$, the matrix of $(F_1^\sharp,F^\sharp)$ gets 
transposed, inverted, and the roles of $L$, $T$ exchanged.

Since $V^\sharp(F_1x)=1\otimes x$ by definition, 
\eqref{EqB2} and \eqref{EqB3} imply that
$F^{\prime\sharp}=(V^\sharp)^\vee$ and $V^{\prime\sharp}=(F^\sharp)^\vee$, 
thus dualising interchanges etale and multiplicative as well as 
$V$-nilpotent and $F$-nilpotent Dieu\-donn\'e\ displays.

\subsection*{Definition of the duality homomorphism.}

In the remainder of this article, 
let always $G=\BT(\PPP)$ and $G'=\BT(\PPP')$ and, 
by a slight abuse of notation, $G^t=\BT(\PPP^t)$.
For arbitrary $\PPP$ and $\PPP'$ we want to construct 
a functorial homomorphism
$$
\psi:\Bil(\PPP'\times\PPP,\GGG_m)\to
\Ext^1(G'\otimes^LG,\hGm).
$$
Given a bilinear form $\alpha:\PPP'\times\PPP\to\GGG_m$
let us first define a map of complexes 
$\gamma:Z(\PPP')\otimes Z(\PPP)\to Z(\GGG_m)$, which is
equivalent to homomorphisms $\gamma_0$ and $\gamma_1$ 
forming the following commutative diagram, 
where $\varphi=F_1-\id$ and $\varphi'=F_1'-\id$. 
$$
\xymatrix@C+3.5em@M+0.2em@R-0em{
Q'\otimes Q \ar[r]^-{\id\otimes \varphi+\varphi'\otimes\id} 
\ar[d]^{\gamma_0} & 
Q'\otimes P\oplus P'\otimes Q 
\ar[r]^-{-\varphi'\otimes\id+\id\otimes \varphi} 
\ar[d]^{\gamma_1} & P'\otimes P \\
{\II_R} \ar[r]^{v^{-1}-\id} & {\WW(R)} 
}
$$
We let $\gamma_0=\alpha$, and there are 
\label{PaCh}
two choices for $\gamma_1$: either  
$\gamma_1(q'\otimes p+p'\otimes q)=\alpha(q',p)+\alpha(p',F_1q)$ or 
$\gamma_1(q'\otimes p+p'\otimes q)=\alpha(F_1'q',p)+\alpha(p',q)$.
The two resulting maps $\gamma$ are homotopic via the homotopy
consisting of $\alpha:P'\otimes P\to\WW(R)$ in top degree and zero 
in lower degrees, in particular the homomorphism
in $D(\tCCC_R)$
$$
\gamma^L:Z(\PPP')\otimes^L Z(\PPP)\xrightarrow{\can} 
Z(\PPP')\otimes Z(\PPP)\xrightarrow{\gamma}Z(\GGG_m)
$$
is independent of the choice. From $\gamma^L$ we obtain 
the following morphism in $D(\tCCC_R)$
that can be viewed as an element of 
$\Ext^1(G'\otimes^LG,\hGm)$; 
by definition this is $\psi(\alpha)$.
\begin{multline*}
G'\otimes^L G\cong
Z(\PPP')\otimes^L Z(\PPP)\otimes^L\QQ_p/\ZZ_p\otimes^L\QQ_p/\ZZ_p \\
\xrightarrow{\gamma^L\otimes\id\otimes\id}
Z(\GGG_m)\otimes^L\QQ_p/\ZZ_p\otimes^L\QQ_p/\ZZ_p\cong
\hGm\otimes^L\QQ_p/\ZZ_p\overset{\eqref{Eq2}}{\cong}\hGm[1]
\end{multline*}

A direct computation of $G'\otimes^L G$ yields the 
following isomorphism, explained in detail in SGA 7, VIII, 1.3.
\begin{equation}
\label{Eq6}
\Ext^1(G'\otimes^LG,\hGm)\cong\Hom(G',G^\vee)
\end{equation}

\begin{Def}
The duality homomorphism $\Psi:G^t\to G^\vee$ 
is the image of the canonical bilinear form 
$\PPP^t\times\PPP\to\GGG_m$ under $\psi$ composed with \eqref{Eq6}.
\end{Def}

\begin{Remark*}
The homomorphism $\Psi$ is compatible with base change,
see section \ref{Se4}.
\end{Remark*}

Naturality of $\psi$ gives the next commutative diagram,
which is just an explication of the fact that $\psi$ and
$\Psi$ are equivalent by the Yoneda lemma.
$$
\xymatrix@M+0.2em@R-0em{
{\Hom(\PPP',\PPP^t)} \ar[r]^-{\BT} \ar@{=}[d] &
{\Hom(G',G^t)} \ar[d]^{\Psi_*} \\
{\Bil(\PPP'\times\PPP,\GGG_m)} \ar[r]^-{\psi} &
{\Hom(G',G^\vee)}
}
$$
Here by Theorem \ref{Thm1} we know (but shall not use) 
that $\BT$ is bijective.

\begin{Lemma}
\label{Le5}
There is the following commutative diagram:
$$
\xymatrix@M+0.2em@R-0em{
{\Bil(\PPP'\times\PPP,\GGG_m)}
\ar[r]^-{\psi} \ar@{-}[d]^{\cong} &
{\Ext^1(G'\otimes^LG,\hGm)}
\ar[r]^-{\eqref{Eq6}} \ar@{-}[d]^{\cong} &
{\Hom(G',G^\vee)} \ar@{-}[d]^{\cong} \\
{\Bil(\PPP\times\PPP',\GGG_m)}
\ar[r]^-{\psi} &
{\Ext^1(G\otimes^LG',\hGm)} \ar@{}[ul]|{-1}
\ar[r]^-{\eqref{Eq6}} &
{\Hom(G,G^{\prime\vee})} \ar@{}[ul]|{-1}
}
$$
\end{Lemma}

\begin{proof}
The right hand square is SGA 7, VIII, Proposition 2.2.11.
The left hand square follows from the definitions, 
the sign resulting from the transposition automorphism of 
$\QQ_p/\ZZ_p\otimes^L\QQ_p/\ZZ_p$ which is $-\id$.
\end{proof}

In particular, a skew symmetric bilinear form 
$\PPP\times\PPP\to\GGG_m$ gives a symmetric biextension
and an anti-symmetric homomorphism $G\to G^\vee$.

\begin{Thm}
\label{Thm2}
For every Dieudonn\'e\ display $\PPP$ over an admissible
local ring $R$ the duality homomorphism 
$\Psi:\BT(\PPP^t)\to\BT(\PPP)^\vee$ is an isomorphism. 
\end{Thm}

By Theorem \ref{Thm1} this implies that $\psi$ is 
an isomorphism as well.

For the proof of Theorem \ref{Thm2} we begin with a
number of reductions.
Using the decompositions \eqref{EqS1} and \eqref{EqS2} we may assume 
that $\PPP$ is etale or of multiplicative type or bi-nilpotent.
By Lemma \ref{Le5} the assertion for $\PPP$ is equivalent 
to the assertion for $\PPP^t$,
so the multiplicative case can be omitted.
Since a homomorphism of $p$-divisible groups 
over $R$ is an isomorphism if and only if it is an
isomorphism over $\bar k$, we may assume that
$R$ is an algebraically closed field.

\begin{proof}[Proof of Theorem \ref{Thm2} in the etale case]
We may assume $G=\QQ_p/\ZZ_p$, accordingly 
$\PPP=(\WW(R),\WW(R),pf,f)$ and $\PPP^t=\GGG_m$.
In order that $\Psi:\hGm\to(\QQ_p/\ZZ_p)^\vee$ is an 
isomorphism it suffices that it is not divisible by $p$.
To get $\Psi$ we trace the definition of $\psi$
applied to the natural bilinear form 
$\alpha:\GGG_m\times\PPP\to\GGG_m$.

Under the quasi-isomorphisms $Z(\GGG_m)\simeq\hGm[-1]$ 
and $T_p(\QQ_p/\ZZ_p)\simeq Z(\PPP)$ the homomorphism 
$\gamma^L$ (defined by the first choice for $\gamma_1$) 
gets identified with the tautological isomorphism 
$\hGm[-1]\otimes^L T_p(\QQ_p/\ZZ_p)\cong\hGm[-1]$. 
Hence $\psi(\alpha)$ is the isomorphism 
$\hGm\otimes^L\QQ_p/\ZZ_p\cong\hGm[1]$
induced by \eqref{Eq2}, in particular $\psi(\alpha)$ is not
divisible by $p$, so the same is true for $\Psi$,
its image under \eqref{Eq6}. 

Note that we did not use the definition of \eqref{Eq6}.
We leave it to the reader to determine whether 
$\Psi$ is the identity or its negative.
\end{proof}

The bi-nilpotent case relies on the following theorem of Cartier.
For any ring $A$ let $\hatW$ be the functor on $A$-algebras
$\hatW(B)=\hatW(\NNN(B))$. Then by \cite{Cartier} the bilinear maps
\begin{equation*}
W(A)\times\hatW(B)\xrightarrow{\mult}
\hatW(B)\xrightarrow{\hex}\hGm(B)
\end{equation*}
induce an isomorphism $W(A)\cong\Hom(\hatW,\hGm)$,
thus an isomorphism of sheaves on $\CCC_R^{\op}$ 
\begin{equation}
\label{Eq7}
W\cong\HHom(\hatW,\hGm)
\end{equation} 
because passing from functors on all $R$-algebras 
to functors on $\CCC_R$ makes no difference as
$\hatW$ is the direct limit of functors 
that are represented by rings in $\CCC_R$.
Let $W[f]$ be the kernel of $f:W\to W$.
Since $f$ corresponds to the dual of $v$ under \eqref{Eq7} 
and the cokernel of $v:\hatW\to\hatW$ is $\hGa$, 
we deduce an isomorphism
\begin{equation}
\label{Eq8}
W[f]\cong\Hom(\hGa,\hGm).
\end{equation}

\begin{Lemma}
\label{Le6}
The Frobenius homomorphism $f:W\to W$ defines a surjective 
homomorphism of flat sheaves on $\CCC_R^{\op}$.
\end{Lemma}

\begin{proof}
As a functor on rings
$W$ is represented by $B=\ZZ[X_0,X_1,\ldots\,]$ 
and $f$ by a faithfully flat ring homomorphism $f^\sharp:B\to B$.
In fact, since the truncated Frobenius $f_n:W_{n+1}\to W_n$ 
is a group homomorphism which is surjective on
geometric points, its fibres are one-dimensional. 
Hence $f_n^\sharp:\ZZ[X_0,\ldots,X_n]\to\ZZ[X_0,\ldots,X_{n+1}]$ 
is faithfully flat by \cite{Matsumura} Theorem 23.1,
so $f^\sharp=\varinjlim f_n^\sharp$ is
faithfully flat.
It remains to show that if $B\to S$ is a ring homomorphism 
with $S\in\CCC_R$, then $B\otimes_{f^\sharp,B}S$ 
lies in $\CCC_R$. 
This can be deduced from the relation
$f^\sharp(X_i)\equiv X_i^p$ modulo $p$.
\end{proof}

\begin{proof}[Proof of Theorem \ref{Thm2} in the bi-nilpotent case]
Assume that $\PPP$ is bi-nil\-potent.
Since $\Psi$ is a homomorphism of $p$-divisible groups of the
same height, it suffices that $\Psi$ is injective, or even that
$\Delta\circ\Psi$ is injective for some homomorphism 
$$
\Delta:G^\vee\to\EExt^1(G,\hGm).
$$

Let $\FFF\PPP^t=(\tP^t,\tQ^t,F^t,F_1^t)$ be the display associated 
to $\PPP^t$ and view $\tP^t,\tQ^t$ as sheaves on $\CCC_R^{\op}$.
From \eqref{Eq7} we get an isomorphism 
$u:\tP^t\cong\HHom(P_\NNN,\hGm)$.
We claim that $\Delta$ can be chosen such that we have the 
following commutative diagram in $\tCCC_R$ with exact rows, 
where $i:Q_\NNN\to P_\NNN$ denotes the inclusion.
The diagram is similar to \cite{Zink-Disp} (223) 
in a different technical context.
\begin{equation}
\label{Eq9}
\xymatrix@M+0.2em@R-0em@C-0em{
0 \ar[r] &
{\tQ^t} \ar[r]^{1-F_1^t} \ar[d]^{u\circ F_1^t} &
{\tP^t} \ar[r]^{\pi} \ar[d]^{i^*\circ u} &
G^t \ar[r] \ar[d]^{\Delta\circ\Psi} & 0 \\
0 \ar[r] &
{\HHom(\PPP_\NNN,\hGm)} \ar[r]^{(F_1-1)^*} &
{\HHom(Q_\NNN,\hGm)} \ar[r]^{\delta} &
{\EExt^1(G,\hGm)}
}
\end{equation}

Before proving the claim, let us apply the snake lemma 
and deduce that $\Delta\circ\Psi$ is injective as required.
Since $u$ is an isomorphism it suffices that 
$F_1^t$ is surjective and $u\circ(1-F_1^t)$,
or equivalently $u$, induces an isomorphism 
$\Ker(F_1^t)\cong\Ker(i^*)$.
In terms of a normal decomposition $P=L\oplus T$
and the induced normal decomposition $\tP^t=\tT^t\oplus\tL^t$,
the homomorphism $F_1^t$ is equal to 
\begin{align*}
I_R\tL^t\oplus\tT^t & \longrightarrow (\tL^t)^{(1)}\oplus (\tT^t)^{(1)}
\xrightarrow{(F^t,F_1^t)^\sharp} \tP^t \\
(v(w_1)l,w_2t) & \longmapsto (w_1\otimes l,f(w_2)\otimes t)
\end{align*}
where $(F^t,F_1^t)^\sharp$ is an isomorphism.
By Lemma \ref{Le6} it follows that $F_1^t$ is surjective,
moreover $\Ker(F_1^t)=W[f]\otimes_WT^t$.
Let $\hatV$ be the formal completion of the vector group
$P/Q\cong T/\II_RT$. Then $\Ker(i^*)\cong\HHom(\hatV,\hGm)$,
and \eqref{Eq8} finishes the proof.

Let us now look at those parts of \eqref{Eq9} 
that do not involve $\Delta$.
It is straightforward that the left hand
square commutes, using that $\hex(v(w))=\hex(w)$
for $w\in\hatW$ according to \eqref{Eqhex}.
Consider the following complexes of sheaves in
$\tCCC_R$ concentrated in degrees $0,1$.
$$
Z=[Q_\NNN\xrightarrow{F_1-1}P_\NNN],\qquad 
Z^t=[\tQ^t\xrightarrow{F_1^t-1}\tP^t]
$$
Since $\PPP^t$ is $V$-nilpotent, 
$Z(\PPP^t)$ is quasi-isomorphic to
$G^t[-1]$ by Proposition \ref{Pr2}.
Since $\PPP^t$ is also $F$-nilpotent,
the inclusion $Z(\PPP^t)\to Z^t$ is a 
quasi-isomorphism by \cite{Zink-Disp} Corollary 82.
This gives the exact upper row of \eqref{Eq9}.
The lower row arises from the short exact sequence 
$$
0\to Q_\NNN\xrightarrow{F_1-\incl}P_\NNN\to G\to 0,
$$
which exists because $\PPP$ is $V$-nilpotent.
Here $\HHom(G,\hGm)$ vanishes because
$\PPP$ is $F$-nilpotent, thus $G$ unipotent;
see Corollary \ref{Cor1}.

Finally, let us define $\Delta$ to be the image
of $\id_{G^\vee}$ under the first row of the following
commutative diagram, whose horizontal arrows $\beta$
are given by adjunction.\footnote{%
A simpler definition of $\Delta$,
equivalent to the above according to
SGA 7, VIII, Proposition 2.3.11 and its 
correction \cite{BBM} p.~253, is the following: 
The restriction of $\Delta$ to $G^\vee[p^r]=G[p^r]^\vee$ 
is the connecting homomorphism of the $\Ext$-sequence
associated to $0\to G[p^n]\to G\to G\to 0$. 
As $G$ is unipotent in our case, it follows that
$\Delta$ is an isomorphism by \cite{Mazur-Messing} Theorem 10.2
or by Proposition \ref{PrX}.}
The rest of the diagram is used to determine
the composition $\Delta\circ\Psi\circ\pi$ and
show that the right hand square of \eqref{Eq9}
commutes.
$$
\xymatrix@M+0.2em{
{\Hom(G^\vee,G^\vee)} \ar[r]^-{\cong}_-{\eqref{Eq6}} \ar[d]_{\Psi^*} &
{\Ext^1(G^\vee\otimes^LG,\hGm)} \ar[r]^-{\beta_1} \ar[d]^{(\Psi\otimes\id)^*} &
{\Hom(G^\vee,\EExt^1(G,\hGm))} \ar[d]^{\Psi^*} \\
{\Hom(G^t,G^\vee)} \ar[r]^-{\cong}_-{\eqref{Eq6}} &
{\Ext^1(G^t\otimes^LG,\hGm)} \ar[r]^-{\beta_2} \ar[d]^{(\pi\otimes\id)^*} &
{\Hom(G^t,\EExt^1(G,\hGm))} \ar[d]^{\pi^*} \\
&
{\Ext^1(\tP^t\otimes^LG,\hGm)} \ar[r]^-{\beta_3} &
{\Hom(\tP^t,\EExt^1(G,\hGm))} 
}
$$
By the definition of $\Psi$, the image of $\id_{G^\vee}$ in the
middle is equal to $\psi(\alpha)$, where $\alpha$ is the
natural bilinear form $\PPP^t\times\PPP\to\GGG_m$.

Because of the quasi-isomorphisms 
$Z\to Z(\PPP)$ and $Z(\PPP^t)\to Z^t$,
in the construction of $\psi(\alpha)$ 
we can start with the obvious pairing 
$\gamma:Z^t\otimes Z\to Z(\hGm)\simeq\hGm[-1]$,
defined by the second choice of $\gamma_1$ 
on page \pageref{PaCh}.
Since the double tensor product $\otimes^L\QQ_p/\ZZ_p$
results in a shift by two, $\psi(\alpha)$ gets 
identified with the composition
$$
G^t\otimes^LG\xleftarrow{\simeq}
Z^t[1]\otimes^L Z[1]\xrightarrow{\can}
Z^t[1]\otimes Z[1]\xrightarrow{\gamma[2]}
\hGm[1].
$$
Using that $\pi$ is induced by the obvious
homomorphism $\tP^t\to Z^t[1]$, it follows that 
$\Delta\circ\Psi\circ\pi=\beta_3((\pi\otimes\id)^*(\psi(\alpha)))$
is equal to the upper line of the following diagram,
where $\gamma'$ is induced by the pairing $\gamma[2]$,
the arrow $\sigma$ is induced by the quasi-isomorphism
$Z[1]\to G$, and $\tau$ is
the obvious homomorphism $Z\to Q_\NNN$.
$$
\xymatrix@M+0.2em{
{\tP^t} \ar[r]^-{\gamma'} \ar[dr]_-{i^*\circ u} &
{\HHom(Z[1],\hGm[1])} \ar[r]^-{\sigma} &
{\EExt^1(G,\hGm)} \\
& {\HHom(Q_\NNN[1],\hGm[1])} \ar[u]_{\tau[1]^*} \ar[ur]_-{\delta}
}
$$
Since both triangles commute, we obtain
$\Delta\circ\Psi\circ\pi=\delta\circ i^*\circ u$
as desired.
\end{proof}

\section{Change of the base ring}
\label{Se4}

Let $f:R\to R'$ be a local homomorphism of local rings
which are admissible in the sense of Definition
\ref{DefAdm}. It is no surprise that all preceding
constructions are compatible with base change by $f$,
i.e.\ for a Dieudonn\'{e} display $\PPP$ over $R$
we have a natural isomorphism
\begin{equation}
\label{EqBaseBT}
u:\BT(\PPP_{R'})\cong\BT(\PPP)_{R'},
\end{equation}
transitive with respect to triples 
$R\to R'\to R''$, such that the following commutes.
\begin{equation}
\label{EqBasePsi}
\xymatrix@M+0.2em{
\BT(\PPP^t_{R'}) \ar[r]^{\Psi} \ar[d]_u &
\BT(\PPP_{R'})^\vee \ar[d]^u \\
\BT(\PPP^t)_{R'} \ar[r]^{\Psi} &
\BT(\PPP)^\vee_{R'}
}
\end{equation}

If one uses the original construction 
of the functor $\BT$ in \cite{Zink-DDisp}, 
the isomorphism $u$ is quite clear, 
but for \eqref{EqBasePsi} we need $u$ in terms
of the formulae of Definition \ref{DefBT}.
This is a question about functoriality
of the category $\CCC_R$.

Assume first that the residue extension of $R\to R'$
is algebraic. Then every $S\in\CCC_{R'}$ lies
in $\CCC_R$ too, and coverings of $S$ in both categories
are the same. Hence we have an exact restriction
functor $\tCCC_{R'}\to\tCCC_R$ and \eqref{EqBaseBT} 
becomes evident. 
By construction of $\psi$ the following diagram
commutes, which gives \eqref{EqBasePsi}.
$$
\xymatrix@C+0.2em{
\Bil(\PPP'\times\PPP,\GGG_m) \ar[r]^-{\psi} \ar[d] &
\Hom(\BT(\PPP')\otimes^L\BT(\PPP),\hGm[1]) \ar[d] \\
\Bil(\PPP'_{R'}\times\PPP_{R'},\GGG_m) \ar[r]^-{\psi} &
\Hom(\BT(\PPP'_{R'})\otimes^L\BT(\PPP_{R'}),\hGm[1])
}
$$

In general we have to modify $\CCC_R$
in order to apply the same reasoning. 
Let $\EEE_R$ be the category of all $R$-algebras
$S$ such that the nilradical $\NNN_S$ is nilpotent,
$\NNN_S$ contains $\Fm S$, and $S_{\red}=S/\NNN_S$
is perfect. Let $\tEEE_R$ be the category of abelian
sheaves on $\EEE_R^{\op}$ for the topology where
a covering is a faithfully flat homomorphism
$S\to S'$ such that $S_{\red}\to S'_{\red}$
is ind-etale. 
The last condition is automatic
when $S$ and $S'$ lie in $\CCC_R$;
conversely for $S\in\CCC_R$ and a covering
$S\to S'$ in $\EEE_R$ we necessarily have
$S'\in\CCC_R$.
It follows that coverings of $S\in\CCC_R$ are
the same in $\CCC_R$ or in $\EEE_R$,
whence an exact restriction functor
$\tEEE_R\to\tCCC_R$.

Now we simply note that in all constructions
we could use $\EEE_R$ in place of $\CCC_R$. 
Then \eqref{EqBaseBT} and \eqref{EqBasePsi}
follow as before since every $S\in\EEE_{R'}$ 
lies in $\EEE_R$
with the same coverings in both categories.
The only point
that might need verification is the fact
that $p$-divisible groups over $R$ form an
{\em exact} subcategory of $\tEEE_R$.
This follows from the remarks preceding
Corollary \ref{CorTop} or from:

\begin{Lemma}
If $H$ is a finite flat group scheme over $R$
and $\Spec T\to\Spec S$ is an $H$-torsor with
$S\in\EEE_R$ then $S\to T$ is a covering in $\EEE_R$.
\end{Lemma}

\begin{proof}
We may assume that $H$ is etale or $H_k$ is 
infinitesimal since any $H$ is an extension of
such groups. The etale case is clear.
If $H_k$ is infinitesimal then $H_S$ is
infinitesimal too, so $T_{\red}\cong S_{\red}$
as $S_{\red}$ is perfect.
\end{proof}

\appendix

\section{Infinitesimal extensions of $p$-divisible groups}

For a lack of reference let us mention the following,
probably well-known, generalisation of the deformational duality
theorem in \cite{Mazur-Messing}.
Suppose $G$, $H$ are $p$-divisible groups on an arbitrary 
scheme $S$ and $S_\circ$ is a closed subscheme of $S$.
Let $\Hom(T_pG,H)=\varinjlim\Hom(G[p^n],H)$
with transition maps induced by $p:G[p^{n+1}]\to G[p^n]$.
We have a homomorphism
$$
\delta:\Hom_{S/S_\circ}(T_pG,H)\to\Ext^1_{S/S_\circ}(G,H)
$$
induced by the exact sequences 
$0\to G[p^n]\to G\xrightarrow{p^n}G\to 0$,
where $\Hom_{S/S_\circ}$ denotes homomorphisms on $S$ which are
trivial on $S_\circ$ and $\Ext^1_{S/S_\circ}$ denotes 
isomorphism classes of extensions on $S$ 
equipped with a trivialisation on $S_\circ$.

\begin{Prop}
\label{PrX}
If the quasicoherent ideal $I\subseteq\OOO_S$ defining 
$S_\circ$ is nilpotent and annihilated by a power of $p$, 
then $\delta$ is bijective.
\end{Prop}

For $H=\mu_{p^{\infty}}$ this results in an isomorphism
$G^\vee(S/S_\circ)\cong\Ext^1_{S/S_\circ}(G,\mu_{p^\infty})$,
which is \cite{Mazur-Messing} Theorem 10.2, 
but the isomorphism given
there is the negative of $\delta$ by Lemma \ref{LeX} below.

\begin{proof}[Proof of Proposition \ref{PrX}]
Assume that $p^rI=0$ and $I^n=0$ and let $m=nr$. 
The inverse of $\delta$ can be constructed as follows.
Assume that $e=[H\to E\to G]$ is an extension 
on $S$ trivialised on $S_\circ$, i.e.\ provided with a section 
$s_\circ:G_\circ\to E_\circ$. 
Then $p^ms_\circ$ lifts to a unique homomorphism
$t:G\to E$, giving the following morphism of exact sequences.
$$
\xymatrix@M+0.2em@R-0em{
0 \ar[r] & 
H \ar[r] \ar[d]^{\id} & 
H\times G \ar[r] \ar[d]^{(\id,t)} &
G \ar[r] \ar[d]^{p^m} & 0 \\
0 \ar[r] & H \ar[r] & E \ar[r] & G \ar[r] & 0
}
$$
The kernel of $(\id,t)$ is the graph of a homomorphism 
$f:G[p^m]\to H$ that is trivial on $S_\circ$, 
and $e\mapsto -f$ is the inverse of $\delta$.
\end{proof}

\begin{Lemma}
\label{LeX}
The natural diagram
$$
\xymatrix@M+0.2em{
\Hom_{S/S_\circ}(T_pG,H) \ar[r]^-\delta \ar@{-}[d]^\cong &
\Ext^1_{S/S_\circ}(G,H) \ar@{-}[d]^\cong \\
\Hom_{S/S_\circ}(T_pH^\vee,G^\vee) \ar[r]^-\delta &
\Ext^1_{S/S_\circ}(H^\vee,G^\vee)
}
$$
whose vertical isomorphisms are given by duality is anti-commutative.
\end{Lemma}

\begin{proof}
Let $f\in\Hom_{S/S_\circ}(T_G,H)$ be given and let $H\to E\to G$ be
its image under $\delta$. For every sufficiently large $n$ so that $f$ 
is represented by a homomorphism $f_n:G[p^n]\to H[p^n]$, the 
truncated extension $E[p^n]$ is naturally isomorphic to the middle
cohomology of the following complex, denoted $K(f_n)$.
$$
G[p^n]
\xrightarrow{(\id,-f_n)}
G[p^n]\oplus H[p^n]
\xrightarrow{(f_n,\id)}
H[p^n]
$$
Since $K(f_n)^\vee\cong K(-f_n^\vee)$ the assertion follows.
\end{proof}

\end{document}